\documentclass{elsartA}

 \usepackage{graphicx}

\usepackage{amssymb}
\newtheorem{definition}{Definition}
\newtheorem{conjecture}{Conjecture}
\newcommand{\R}{{\Bbb R}}

\newcommand{\N}{{\Bbb N}}
\newcommand{\C}{{\Bbb C}}

\begin{document}
\begin{frontmatter}
\title{Slowly oscillating wavefronts of the KPP-Fisher delayed equation}
\author[a]{Karel Hasik}
\author[a,c]{and Sergei Trofimchuk}
\address[a]{Mathematical Institute, Silesian University, 746 01 Opava, Czech Republic}
\address[c]{Instituto de Matem\'atica y Fisica, Universidad de Talca, Casilla 747,
Talca, Chile \\ {\rm Karel.Hasik@math.slu.cz and
trofimch@inst-mat.utalca.cl}}

\bigskip

\begin{abstract}
\noindent  This paper concerns the semi-wavefronts (i.e. bounded  solutions  $u=\phi(x\cdot\nu +ct) >0,$ $ |\nu|=1, $ satisfying $\phi(-\infty)=0$) 
to the delayed KPP-Fisher equation
$$u_t(t,x) = \Delta u(t,x)  + u(t,x)(1-u(t-\tau,x)), \ u \geq 0,\ x
\in \R^m. \eqno(*)$$  First, we show that each semi-wavefront should be either monotone or slowly oscillating.  Then a complete solution to the problem of existence of  semi-wavefronts is provided.  We prove next that the semi-wavefronts are in fact wavefronts (i.e. additionally $\phi(+\infty)=1$) if $c \geq 2$ and $\tau \leq 1$; our proof  uses dynamical properties of some auxiliary one-dimensional  map with the negative Schwarzian. 
The analysis of the fronts' asymptotic expansions at infinity is another key ingredient of our approach. It allows to indicate the maximal domain ${\mathcal D}_n$ of  $(\tau,c)$ where the existence of non-monotone wavefronts can be expected.  Here we show that the problem of  
wavefront's existence  is closely related to the  Wright's  global stability conjecture. 

\end{abstract}
\end{frontmatter}
\newpage

\section{Introduction and main results}\vspace{-3mm}
\label{intro} 

\vspace{-3mm} The delayed KPP-Fisher equation or the
diffusive Hutchinson's equation
\begin{equation}\label{17} \hspace{5mm}
u_t(t,x) = \Delta u(t,x)  + u(t,x)(1-u(t-\tau,x)), \ u \geq 0,\ x \in
\R^m, \end{equation} 

\vspace{-6mm} 

can be considered as one of the most important 
examples of delayed reaction-diffusion equations. In particular, during the past decade, this model  
has been studied by many authors, see 
\cite{ZAMP,FW,fhw,FGT,FTnl,GT,KO,wz} and the references therein.  A significant 
part of the research dealt  with  the existence of  traveling fronts 
connecting the trivial and positive
steady states in (\ref{17}) and in its non-local variant  \cite{BNPR,FZ,GL,VP}

\vspace{-4mm} 

\begin{equation}\label{17nl} \hspace{-7mm}
u_t(t,x) = \Delta u(t,x)  + u(t,x)\left(1-\int_\R K(y) u(t,x-y)dy\right), \ \int_\R K(s)ds=1. 
\end{equation}

\vspace{-5mm}

We recall that the classical solution $u(x,t) = \phi(\nu \cdot  x +ct),$ $|\nu| =1,$ is a
wavefront (or a traveling front) for (\ref{17}) or (\ref{17nl}) propagating at the velocity $c\geq 0$, if the profile
$\phi$ is non-negative and satisfies $\phi(-\infty) = 0$ and
$\phi(+\infty) = 1$. By replacing condition $\phi(+\infty) = 1$ with less restrictive  
$0< \liminf_{s \to +\infty}\phi(s) \leq \limsup_{s \to +\infty}\phi(s) < \infty$, we get the definition of a semi-wavefront. 
The non-negativity  requirement $\phi \geq 0$  is due to  the biological interpretation of $u$ as of the concentration of a  dominant gene  that is reminiscent of the  seminal works by Kolmogorov, 
Petrovskii, Piskunov and Fisher. 

\vspace{-1mm}

Recently, the wavefront existence problem  for (\ref{17}), (\ref{17nl})  was
considered by using quite different approaches. The first method
was proposed by Wu and Zou in \cite{wz}. It uses the positivity and
monotonicity properties of the integral operator
\vspace{-5mm} 

\begin{equation}\label{psea} \hspace{-7mm}
(A\phi)(t) = \frac{1}{z_2-z_1}\left\{\int_{-\infty}^te^{z_1
(t-s)}(\mathcal{H}\phi)(s)ds + \int_t^{+\infty}e^{z_2
(t-s)}(\mathcal{H}\phi)(s)ds \right\},
\end{equation}
\vspace{-5mm} where $(\mathcal{H}\phi)(s)= \phi(s)(b+1
-\phi(s-h)), \ h := c\tau,$   is taken with some appropriate $b >1$, 
and  $z_1<0<z_2$ satisfy $z^2 -cz
-b =0$. A direct verification shows that the
profiles $\phi \in C(\R,\R_+)$ of semi-wavefronts can be also identified as positive bounded solutions of  the integral equation $A\phi=\phi$ satisfying 
the above mentioned boundary conditions at $\pm\infty$.  Unfortunately, the presence of 
positive delay in (\ref{psea}) strongly affects the monotonicity of $A$. In order to overcome this difficulty, 
two different orderings, the usual  one  and a non-standard Smith and Thieme ordering of  $C(\R,\R_+)$, were combined in  \cite{wz}. 
Even so  the operator $A$ was monotone with respect to each of  these two orderings only for sufficiently 
small $h$  and  monotone $\phi$.  

\vspace{-1mm}

The operator $A$ is well defined when $b>0$. Taking formally $b=-1$ in (\ref{psea}) and interpreting  
correctly the obtained expression for $c>2$, \mbox{instead of $A$ we obtain }
\begin{equation}\label{iie} \hspace{5mm}
(B\varphi)(t) = \frac{1}{\mu-\lambda}
\int_t^{+\infty}(e^{\lambda (t-s)}- e^{\mu
(t-s)})\varphi(s)\varphi(s-h)ds,
\end{equation}

\vspace{-10mm}

where $0 < \lambda < \mu$ are the roots of $z^2 -cz + 1 =0$.   Remarkably,  all monotone
wavefronts to equation (\ref{17}) can be found via a monotone iterative algorithm
which uses $B$ (or its limit version $B_2$ if $c=2$) and converges uniformly on $\R$, see \cite{GT}.  
Similar ideas  were also  successfully applied in \cite{FW,FZ,KO}.  However, our attempts to 
use   the monotone operator $B$ in the case of non-monotone waves were not fruitful. 

Aiming  to get rid of  monotonicity requirements,  Shiwang Ma achieved 
an important progress in  \cite{ma,ma1}. He showed 
that operators similar to $A, B$ have  good compactness properties in suitable Banach spaces. Therefore, in certain situations,   the Schauder fixed point theorem  could be used instead of the iterative monotone scheme from \cite{wz}. Ma's idea was successfully applied to various  reaction-diffusion models with bounded nonlinearities.  Nevertheless, equation  $A\phi=\phi$  with $A$ defined by (\ref{psea}) has never been considered within the Ma's approach: this is mainly   
because of the considerable difficulties related to the construction of a nontrivial $A$-invariant set suitable for the application of the Schauder fixed point theorem.

It is therefore tempting, in order to avoid the construction of a non-trivial bounded  $A$-invariant convex closed set $\Omega$,  to apply the Leray-Shauder continuation principle to equation $A\phi =\phi$.  The main obstacle for the realization of such an idea is the apparent  impossibility to have at the same time  complete continuity of $A$ and the non-empty interior of $\Omega$.  This problem was avoided in a  nice way by Berestycki {\it et al.}  in \cite{BNPR}. Working with 
equation (\ref{17nl}), for a fixed $\delta >0$,  Berestycki {\it et al.}  considered a family of associated boundary value problems, with the boundary conditions $\phi_n(-n)=0, \phi_n(n)=1, \phi_n(0)= \delta$.  Fortunately,  the above mentioned contradiction between the compactness of operator and the openness of its domain  does not occur on finite intervals $[-n,n]$.   Hence, the Leray-Shauder continuation principle (with corresponding calculation of  {\it a priori} {estimates}, degrees etc) can be applied for each \ $n \in \N$. Finally,  the wave profile $\phi$ was obtained in  \cite{BNPR} as the
limit of $\phi_n$ .  The proof of the existence in  \cite{BNPR} is rather technical and non-trivial.  Regrettably, the conditions of  $C^1$-smoothness  of kernel $K$ and especially the positivity of $K(0) >0$  do not allow  use the existence theorem from 
\cite{BNPR} to deduce a similar result for equation (\ref{17}). Indeed, if we take some $\delta-$like sequence of 
kernels $K^{(j)}(s) \to \delta(s-h)$ then  the corresponding sequence of traveling waves $\phi^{(j)}(s)$ could be eventually unbounded in view of {\it a priori} estimates obtained in \cite{BNPR}. 

Our short description of analytical  tools used to prove the wave existence in (\ref{17}), (\ref{17nl}) would be incomplete without mentioning the  Lin-Hale approach to heteroclinic solutions  developed in \cite{fhw, FTnl}.  This method allowed to obtain almost optimal existence results (i.e. $\tau \leq 3/2$ and $c \geq c'$, for some indefinite and large $c'$: see also Fig. 1 and Conjecture \ref{Coj} below) for rapidly traveling fronts.   Nevertheless,  the most interesting in applications  critical waves  were excluded in  \cite{fhw, FTnl}.  Surprisingly, as  the recent work \cite{FGT} shows,  the Lin-Hale method still can be extended to give a complete solution to the problem of  existence of  monotone fronts in several models (including (\ref{17})).   However, the monotonicity of waves is one of crucial assumptions in \cite{FGT} and, at this moment, it is not clear whether it can be dropped.  

After analyzing the above approaches to the existence problem and motivated by \cite{BNPR,ma,wz}, we decided  to work with the equation $A\phi =\phi$. As a result, we elaborated  a 
framework suitable for the application of the Schauder fixed point theorem for an appropriately  modified version of the operator $A$.  Before stating the corresponding existence theorem, let us define several subsets of parameters $(\tau, c) \in \R_+^2$ (see also Figure 1 below):  

\vspace{-7mm}

$$
\mathfrak{D}_s = \{(\tau,c) \in \R_+:  \mbox{there exists a semi-wavefront to Eq. (\ref{17}})\},
$$
$$
\mathfrak{D}_{m} = \{(\tau,c) \in \R_+:  \mbox{there exists a monotone wavefront to Eq. (\ref{17}})\}, 
$$
$$
\mathfrak{D}_{n} = \{(\tau,c) \in \R_+:  \mbox{there exists a non-monotone wavefront to Eq. (\ref{17}})\}.
$$

\begin{thm}[\small Existence criterion for semi-wavefronts] \label{Te1} $\mathfrak{D}_{s} = \{(\tau,c)\in \R^2_+: c \geq 2 \}$.  Furthermore, 
there exist continuous functions $\delta_\pm: \R^2_+ \to (0, +\infty)$ such that 
$\delta_-(\tau,c) < \phi(t) <  \delta_+(\tau,c), \ t \geq Q_0,\ \phi(t) < 1,\ t < Q_0,$ for each semi-wavefront $u=\phi(x\cdot \nu+ct), |\nu| =1,$ and some appropriate $Q_0 = Q_0(\phi)$.  \end{thm}

The proof of Theorem \ref{Te1} requires a detailed study of oscillation/monotonicity  properties of  semi-wavefront 
profiles.  Here we were inspired by geometrical descriptions from  \cite{TT}  of semi-wavefront profiles to the   
Mackey-Glass type delayed reaction-diffusion equation
\begin{equation}\label{17MG}
u_t(t,x) = \Delta u(t,x)  - u(t,x) + g(u(t-\tau,x)), \ u(t,x) \geq
0,\ x \in \R^m.
\end{equation}


It is  known that in the ordinary case (i.e. when $u=u(t)$) models  (\ref{17}), (\ref{17MG}) can be 
considered  within the same family of differential equations governed by linear 
friction (possibly, degenerate) and negative delayed feedback.  Inclusion of the diffusive terms, however,  
makes the similarity between    (\ref{17}) and  (\ref{17MG})  much less  direct. Nevertheless,  it is still possible to prove that the semi-wavefront profiles to (\ref{17})  share all geometric properties established in the case of  equation  (\ref{17MG}). Amazingly, the statements of corresponding assertions become  even sharper while their  proof simplifies: cf. Theorem \ref{T2}, \ref{T3a} below with Theorems 1,3 in \cite{TT}. 

\begin{thm}[\small Monotonicity of the leading edge of semi-wavefronts] \label{T2} \hfill \\
 Let $u(x,t) = \phi(\nu \cdot x +ct),\  |\nu| =1,$ be a non-negavite non constant (possibly, unbounded)  
solution of equation (\ref{17}) satisfying $\phi(-\infty)=0$. Then $\phi(t) > 0, $ $ t \in\R,$ and $\phi$ has a monotone leading edge. The latter means that $\phi'(s) > 0$ on $(-\infty, T_0)\cup (T_1, T_2)$
and $\phi'(s) < 0$ on $(T_0,T_1)$ for some $T_2 \geq T_1 \geq T_0 \in \R \cup \{+\infty\}$. Furthermore, $T_0$ is
finite if and only if $\phi(T_0) > 1$. 
\end{thm}
Before stating the next theorem, we need to introduce the concepts  of critical speeds $c^* < c^\star$  and 
slowly oscillating semi-wavefronts for equation (\ref{17}). Let   $\psi(z,c):=
z^2-cz-\exp(-z c\tau)$ and set $\tau_1:= 0.560771160\dots$

\begin{figure}[t]\label{F1}
\begin{center}
\includegraphics[scale=0.7]{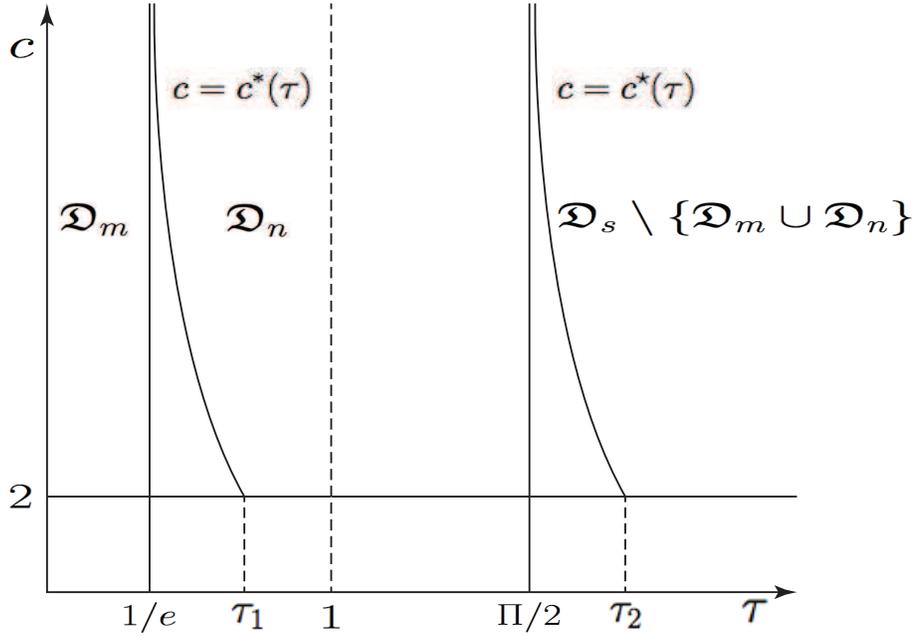}
\end{center}
\caption{Schematic presentation of the critical speeds and delays.}
\end{figure}

By \cite[Lemma 3]{GT},  there exists function $c^*=c^*(\cdot): [0, \tau_1] \to [2,+\infty]$ such that 
$\psi(z,c), \ c \geq 2,$  has has exactly two (counting multiplicity) negative zeros 
in the half plane $\{\Re z < 0\}$ if and only if $\tau \in [0, \tau_1]$ and $c \in [2, c^*(\tau)]$.  Moreover, 
$c^*(\tau) = +\infty$ if and only if $\tau \leq 1/e$ while on the interval $(1/e,\tau_1]$ function $c^*(\tau)$ is decreasing and $c^*(\tau_1) =2$, see Fig. 1.  Similarly, we have the following 
\begin{lem}\label{c1de}  Let $c \geq 2, \ \tau \geq 0$. Then $\psi(z,c)$ has exactly one (counting multiplicity) 
zero in the right half-plane $\Re z > 0$ if and only if one of the
following conditions holds
\begin{enumerate}
  \item $0\leq  \tau \leq \pi/2$, and $c \leq c^\star(\tau):=+\infty,$ 
  \item $\pi/2 < \tau \leq  1.86173\dots:=\tau_2$ and $c \leq c^\star(\tau)$, where $c^\star$ is given implicitly by
\end{enumerate}
\begin{equation}\label{fot}
\tau = \frac{\arccos (-w^2(c^\star)) }{c^\star w(c^\star)}, \quad w^2(c)= \frac{\sqrt{c^4+4}-c^2}{2}.
\end{equation}


Furthermore, if $c > c^\star(\tau)$ and 
$\psi(\lambda_j,c)=0$, $\Re\lambda_j \leq 0$, then  $|\Im \lambda_j | > 2\pi/(c\tau)$. 
\end{lem}
As in \cite{TT},  we follow closely  the definition of slow oscillations from  \cite{mps,mps2}:
\begin{definition} Set $h:= c\tau, \mathbb{K} = [-h,0] \cup \{1\}$. For any
$v \in C(\mathbb{K})\setminus\{0\}$ we define the number of sign
changes by $$\hspace{-1mm} {\rm sc}(v) = \sup\{k \geq 1:{\rm  there
\ are \ } t_0 <
 \dots < t_k  \ {\rm such \ that\ }
v(t_{i-1})v(t_{i}) <0 {\rm \ for \ }  i\geq 1\}. $$ We set ${\rm
sc}(v) =0$ if $v(s) \geq 0$ or  $v(s) \leq 0$ for $s \in
\mathbb{K}$. If $\varphi$ is a non-monotone semi-wavefront profile  to (\ref{17}), we set \ $(\bar \varphi_t)(s) =
\varphi(t+s)- 1$ if $s \in [-h,0]$, and $(\bar \varphi_t)(1) =
\varphi'(t)$. We will say that $\varphi(t)$  is sine-like slowly oscillating 
 if graph of $\varphi$ oscillates around $1$ and has exactly one critical point between each two 
 consecutive intersections with the level $1$, and, in addition, for each $t \geq T_0$ ($T_0$ was defined in Theorem \ref{T2}), it holds  that either sc$(\bar
\varphi_t)=1$ or sc$(\bar \varphi_t)=2$.
\end{definition}
Our next result is similar to \cite[Theorem 3]{TT}. In fact,  it is even stronger, since it excludes non-monotone but 
eventually monotone wavefronts to equation (\ref{17}).  As the numerical simulations of \cite[Figure 1]{BNPR} show, this irregular behavior can occur in simple non-local KPP-Fisher equations. We also believe that such kind of  irregular non-monotone wavefronts can be found  in equation (\ref{17MG}).

\begin{thm}[\small Semi-wavefronts are either monotone or slowly
oscillating] 
\label{T3a}  Let $u= \phi(\nu \cdot x +ct)$ be as in Theorem \ref{T2}. Then one of the next options holds

\vspace{0mm}

\begin{enumerate}
\item $\phi$ is monotonically converging to $1$;
\item $\phi$ is sine-like slowly oscillating around $1$  on a finite maximal interval  and, for some $A>0,$ $ t_0$, it holds  $\phi'(t) >0, $ $ \phi(t) > Ae^{ct},$ $ t \geq t_0$;
\item $\phi$ is sine-like slowly oscillating around $1$ and it is bounded. 
\end{enumerate}
\end{thm}
\begin{rem} \label{MPss} By Theorem \ref{T1} below,  each bounded profile $\phi$ has to develop
non-decaying slow oscillations around $1$ for each  $c
> c^\star(\tau)$ and then, due to \cite{mps2},  these oscillations should be  asymptotically sine-like periodic. 
\end{rem}
The final part of this section concerns the determination of  domain $\mathfrak{D}_{n}\subset \R^2_+$.  We recall that $\mathfrak{D}_{s}$ was already found  in Theorem \ref{Te1} while  the complete description of $\mathfrak{D}_{m}$ was given in  \cite{GT}: 
\begin{prop} 
\label{main} $\mathfrak{D}_{m} = \{(\tau,c)\in \R^2_+: 2 \leq c \leq c^*(\tau)\}$. 
Furthermore,  for some appropriate
$\phi_-$ (given explicitly), we have that $\phi = \lim_{j \to
+\infty} B^j\phi_-$ (if $c>2$), and $\phi = \lim_{j \to
+\infty} B_2^j\phi_-$ (if $c=2$), where the convergence is
monotone and uniform on $\R$. Finally, for each fixed
$c\not=c^*(\tau)$, $\phi(t)$ is the only possible monotone profile (modulo
translation).
\end{prop}
As it was recently demonstrated by Fang and Wu in \cite[Theorem 6.2]{FW}, condition  $c\not=c^*(\tau)$ of Proposition \ref{main} can be dropped.  In any case,  the uniqueness in \cite{FW,GT} was established  only within the class of monotone fronts (see also \cite{FZ} for a similar assertion concerning (\ref{17nl})).   Here, by combining the Berestycki-Nirenberg sliding argument \cite{BN} with the approach of \cite{GT},  we obtain the following 

\vspace{1mm}

\begin{thm} \label{ama} Suppose that $(\tau,c) \in \mathfrak{D}_{m}$  and 
$u=\phi_1,\phi_2$ are wavefronts to (\ref{17}). Then $\phi_1(t)\equiv \phi_2(t+\alpha)$ for some $\alpha \in \R$ and $\phi_j'(t) >0, \ t \in \R$. 
\end{thm}
The sliding solutions method  does not work when $(c,\tau) \not\in \mathfrak{D}_{m}$.   However, as the recent works \cite{AGT,FTnl} have showed,   
the uniqueness (up to translation) of the semi-waveronts  to (\ref{17}) is very likely to be true for 
 large speeds. We believe that for each fixed pair $(\tau,c)$ the semi-wavefront solution to equation 
(\ref{17}) is unique (up to a translation) whenever it exists.  

Theorem \ref{ama} is instrumental in proving 
Theorem \ref{T3a} and, whence, in  establishing our last two results:
\begin{thm}[\small Existence of  non-monotone wavefronts]  \label{ENMTW}
$$\mathfrak{D}_{n}\cup \mathfrak{D}_{m}  \supset \mathcal{D}=\{(\tau,c)\in \R^2_+:  0 \leq \tau \leq 1, \ c \geq 2\}. $$  
Moreover if $(\tau,c) \in \mathcal{D}$ then necessarily $\phi(+\infty)=1$.  Hence, for each $\tau\leq 1$  equation (\ref{17})
has at least one semi-wavefront which necessarily is a wavefront. 
\end{thm}
\begin{cor}[\small Absolute uniqueness of  monotone wavefronts] Suppose that $(\tau,c) \in \mathfrak{D}_{m}$  and 
$u=\phi_1,\phi_2$ are semi-wavefronts to (\ref{17}). Then $\phi_1(t)\equiv \phi_2(t+\alpha)$ for some $\alpha \in \R$ and $\phi_j'(t) >0, \ t \in \R$. 
\end{cor}
\begin{thm}[\small Admissible wavefront speeds and non-existence of fronts] \label{T1}\hfill 
\vspace{0mm}
Eq. (\ref{17}) does not have any travelling front (neither monotone
nor non-monotone) propagating at velocity $c > c^\star(\tau)$ or $c < 2$.
\end{thm}
It can bee seen from Proposition \ref{main} and Theorem \ref{T1}  that  $\mathfrak{D}_{n}\subset  \{(\tau,c) \in \R_+^2: c^*(\tau) < c \leq c^\star(\tau)\}.$ Moreover, by Theorem \ref{ENMTW} and  \cite[Theorem 5.1]{FTnl}, 
$\mathfrak{D}_{n}\cup \mathfrak{D}_{m}\supset  [0,3/2]  \times [c', +\infty) \cup [0,1]  \times  [2,+\infty)$ for 
some large $c'$.  See also Figure 1.  In this way, considerations of the present work suggest 
the following  natural criterion for the existence of non-monotone wavefronts in
(\ref{17}):
\begin{conjecture} \label{Coj}
$\mathfrak{D}_{n}= \{(\tau,c) \in \R_+^2: c^*(\tau) < c \leq c^\star(\tau)\}. $
\end{conjecture}
It can be regarded as an extension of the famous Wright's global stability conjecture \cite{TK,ltt}.
Therefore, in our opinion,  it would be very interesting (and, perhaps, very difficult) 
to prove it.  In particular, in the limit case $c=+\infty$,  Conjecture \ref{Coj} is true if the Wright's conjecture is true. 
An important partial result in proving Conjecture \ref{Coj}  would be the following analog of  the Wright's $3/2$-global stability theorem:  $\mathfrak{D}_{n}\cup \mathfrak{D}_{m}\supset  [0,3/2]\times[2, +\infty).$
 
The structure of the remainder of
this paper is as follows. Section \ref{AU} contains the proof of Theorem \ref{ama}. In the third section, we 
describe the geometrical form of semi-wavefronts.  Theorems \ref{ENMTW}, 
\ref{Te1} and \ref{T1} are proved in Sections  \ref{ETW}, \ref{ESW}, \ref{LaT}
respectively. In Appendix, the characteristic function
of the variational equation at the positive steady state is
analyzed. 

\vspace{-5mm}

\section{Absolute uniqueness of monotone wavefronts} \label{AU}

\vspace{-5mm}

Take some $(\tau,c) \in \mathfrak{D}_{m}$. Then by Proposition  \ref{main} and  \cite[Theorem 6.2]{FW} there exists a {\it unique} monotone 
wavefront $u=\psi_2(\nu\cdot  x+ct)$. Suppose that $u=\psi_1(\nu\cdot x+ct)$ is a different (and therefore non-monotone) wavefront.  Clearly, each profile $\psi_i(t)$ satisfies \vspace{-3mm}
\begin{eqnarray} \label{twe2a} &&
\phi''(t) - c\phi'(t) + \phi(t)(1-
\phi(t-h)) =0,  \ h := c\tau,
\\ \nonumber
&& \phi(-\infty) =0,  \quad \phi(t) \geq 0, \quad t \in \R, \end{eqnarray}
and therefore it is strongly positive due to
\begin{lem} \label{po} Let non-negative $\phi \not\equiv 0$ solve (\ref{twe2a}). Then $\phi(t) >0, \ t \in \R$.   
\end{lem}

\vspace{-5mm}

\begin{pf} Suppose that, for some $s$, we have $\phi(s)=0$. Since $\phi(t) \geq 0, \ t \in \R,$ this 
yields $\phi'(s)=0$. Therefore $y=\phi(t)$ satisfies the following initial value problem 
for a linear second order ordinary differential equation
$$
y''(t) -cy'(t) +(1-\phi(t-h))y(t)=0, \quad y(s)= y'(s) =0. 
$$ 
But then $y(t)  \equiv 0$ due to the uniqueness theorem. \hfill $\square$
\end{pf}
We also will need the asymptotical description of profiles $\psi_i$ at $\pm\infty$. Recall that  $0 < \lambda \leq 
\mu$ denote  the roots of $z^2 -cz + 1 =0$, $c \geq 2$. 
\begin{lem} \label{af1}  Let $c >2, \ q \in \R$. Then,  for sufficiently small $\epsilon >0$, it holds 
$$
\psi_i(t+q) = \frac{e^{\lambda (t+q)}}{\sqrt{c^2-4}}\int_{\R}e^{-\lambda s}\psi_i(s)\psi_i(s-h)ds + O(e^{(\lambda+\epsilon)t}), \quad  t \to -\infty. 
$$
Similarly, if $c=2$ then 
$$
\psi_i(t+q) = e^{\lambda (t+q)}\int_{\R}e^{-\lambda s}\psi_i(s)\psi_i(s-h)ds(-t +O(1)), \quad  t \to -\infty. 
$$
\end{lem}

\vspace{-8mm}

\begin{pf} It is a straightforward consequence of \cite[Lemma 28]{GT}. See also proof of Theorem 6 in \cite{GT}. \hfill $\square$
\end{pf}
\begin{lem} \label{af2}  Suppose that $(\tau,c) \in \mathfrak{D}_{m}$ and let $\lambda_0 <0$ be as in Lemma \ref{roots}.  Let $c \in [2,c^*(\tau)), \ q \in \R,$ and $\epsilon >0$ be sufficiently small. Then 
$$
\psi_i(t+q) = 1- K_ie^{\lambda_0 (t+q)} + O(e^{(\lambda_0- \epsilon)t}), \quad  t \to +\infty,
$$
for some $K_2>0$ and $K_1 \in \R$ independent on $q$. Similarly, if $c=c^*(\tau)$ then 
$$
\psi_i(t+q) = 1 - e^{\lambda_0 (t+q)}(K_it + O(1)), \quad  t \to +\infty, \quad  K_2 >0, \ K_1 \in \R.   
$$
\end{lem}

\vspace{-9mm}

\begin{pf} In the monotone case (i.e. $i =2$), this statement follows from \cite[Lemma 28]{GT} and Lemma \ref{roots} (see also \cite[Theorem 6]{GT} for more details). 
Next,  due to  \cite[Lemma 10]{GT},  the condition 
$(\tau,c) \in \mathfrak{D}_{m}$ implies the hyperbolicity of the positive equilibrium  of (\ref{twe2a}) and therefore $|\psi_1(t)-1|$ converges exponentially to $0$ at $+\infty$. With this observation,  
the analysis of the non-monotone wavefront is completely analogous to the monotone case considered in \cite[Section 7]{GT}.  The unique exception is the sign of $K_1$. Indeed,  in virtue of non-monotonicity of the wavefront $\psi_1$,  $K_1$ could take any real value including $0$.   
 \hfill $\square$
\end{pf}

\vspace{-7mm}

By applying a sliding argument,  we are ready now to prove Theorem \ref{ama}. Set
$$
\mathcal{Q} := \{q: \psi_1(t) \geq \psi_2(t-q), \ t \in \R\}. 
$$
It follows from Lemmas  \ref{af1}, \ref{af2} that $\mathcal{Q} \not = \emptyset$.  On the other hand, it is obvious that the set  
$\mathcal{Q} $ is closed, below bounded and connected (the latter is due to the monotonicity of $\psi_2$).  Let 
$q_* = \inf \mathcal{Q}$, we claim that,  for some finite $t_*$, 
\begin{equation}\label{rav}
\psi_1(t_*) = \psi_2(t_*-q_*). 
\end{equation}

\vspace{-5mm}

Indeed, otherwise 
\begin{equation}\label{13}
\psi_1(t) > \psi_3(t):= \psi_2(t-q_*), \quad t \in \R, 
\end{equation}

\vspace{-5mm}

and therefore Lemma \ref{af1} (taken with $q=0$ and applied to $\psi_1$ and $\psi_3$) implies that there are 
$S_0$ and $\delta_0>0$ such that 
$\psi_1(t) > \psi_3(t+\delta),  \ t \leq S_0$ for all  $\delta \in [0, \delta_0]$.  Now, applying Lemma \ref{af2} (with $q=0$) to the profiles $\psi_3$ and $\psi_1$ we obtain that necessarily $K_2 \geq K_1$.  We claim that  $K_2 >K_1$.  Indeed, otherwise $K_2=K_1 >0$ and therefore the uniqueness proof of \cite[Section 6.3]{GT} can be repeated 
for $c \leq c^*(\tau)$, see also  \cite[Theorem 6.2]{FW}.  Hence, $K_2 > K_1$ and therefore there exist $S_1 \geq S_0, \delta_1 >0$ such that $\psi_1(t) > \psi_3(t+\delta),  \ t \geq S_1$ for all  $\delta \in [0, \delta_1]$.   Finally, considering inequality (\ref{13}) on a fixed interval $[S_0,S_1]$, we find that, for some $\delta_2>0$, it holds  
 $\psi_1(t) > \psi_3(t+\delta),  \ t \in [S_0, S_1]$ for all  $\delta \in [0, \delta_2]$. But then 
 $$
 \psi_1(t) > \psi_3(t+\delta),  \ t \in \R,\ \mbox{for all}  \ \delta \in [0, \delta_*], \ \delta_* = \min\{\delta_j, \ j =0,1,2\}.$$
 Therefore $q_* - \delta_* \in \mathcal{Q},$ a contradiction. 
  
  Hence, (\ref{rav}) holds and therefore non-negative function 
  $
  \theta(t) = \psi_1(t)-\psi_3(t)
  $
 attains its zero minimum at $t_*$.  Moreover, as $\theta(t) >0$ for $t \leq S_0$, we may assume that  $t_*$ is the leftmost zero minimum of $\theta$.  Now, it is easy to see that bounded $\theta$ also satisfies the differential 
 equation 
 $$
 \theta''(t) -c \theta'(t) +\theta(t) = \theta(t)\psi_1(t-h) + \theta(t-h) \psi_3(t)=: \Theta(t),  
 $$
so that either
$$
\theta(t) = \frac{1}{\mu-\lambda}
\int_t^{+\infty}(e^{\lambda (t-s)}- e^{\mu
(t-s)})\Theta(s)ds, \quad \mbox{if $c >2$, cf. (\ref{iie})}, $$
$$\  \mbox{or} \quad 
\theta(t) = 
\int_t^{+\infty}(s-t)e^{(t-s)}\Theta(s) ds, \quad \mbox{if $c =2$.} 
$$
Considering the above relations with $t=t_*$, we deduce immediately that $\Theta(s)\equiv 0$ on $[t_*, +\infty)$. However, this can not happen because of  the inequality 
$\Theta(s) \geq \theta(s-h)\psi_3(s) >0, \ s \in [t_*, t_*+h).$ The obtained contradiction ends the proof of Theorem \ref{ama}. 
\hfill $\square$

\vspace{-7mm}

\section{Semi-wavefront's shape}

\vspace{-6mm}

This section contains a detailed analysis of the oscillation and monotonicity properties of profiles $\phi$ corresponding to 
non constant non-negative solutions  $u(t,x)= \phi(\nu\cdot x+ct), |\nu|=1, \phi(-\infty)=0, $  of the delayed KPP-Fisher equation.  The main conclusions of the section (see also Lemma \ref{ogran} below) are presented  as Theorem \ref{T2} and Theorem \ref{T3a}  in Introduction.

By Lemma \ref{po}, similarly to the case of the Hutchinson's equation, the change of variables $\phi(t) = e^{-x(t)}$ can be applied 
to (\ref{twe2a}).  The obtained equation (see equation (\ref{et}) below) is a unidirectional monotone cyclic feedback system with delay \cite{mps2}. Therefore, analogously as it was done in \cite{TT},  fundamental results from \cite{mps,mps2} can be used to demonstrate slowly oscillating character of the non-monotone semi-wavefronts. 
Nevertheless, here we have preferred to give  short and self-contained  direct proofs of this fact,  additionally establishing  sinusoidal  shape of all (and not only periodic as in \cite{mps2}) oscillating solutions.  See also Remark \ref{MPss} in the introduction. 

\begin{lem} \label{o} Let $Q_0$ be such that $0< \phi(s) <1$ for $s < Q_0$ and 
$\phi(Q_0) =1$. Then $\phi'(s) >0$ for all $s \in (-\infty, Q_0]$.  
\end{lem}

\vspace{-7mm}

\begin{pf} If $\phi'(s)=0, \phi(s) \leq 1,$ for some $s \leq Q_0$, then 
necessarily $\phi''(s)<0$ so that $s$ is a critical
point (local maximum) of $\phi$. As a consequence, $\phi'(t) < 0$ for all
 $t> s$  since otherwise there exists $s_1>s$ such that  $\phi'(s_1)<0,$ $\phi''(s_1)= 0,$  $\phi(s_1), \phi(s_1-h) \in (0,1]$, a contradiction. 
 However, $\phi'(t) <0, $ \ $\phi(t) <1,$ $t > s,$ yields $\phi''(t) <0$ for all $t >s$, 
 and therefore $\phi(t)$ can not be positive for large $t$, a contradiction. \hfill $\square$
\end{pf}

\vspace{-7mm}

\begin{lem} \label{mp} Let $Q_0$ be as in Lemma \ref{o} and $Q_1$ be such that $\phi(s) >1$ for all $s$ from 
some maximal open interval $(Q_0, Q_1)$. Then the only options for the geometrical 
shape of $\phi$ on $(Q_0, Q_1)$ are: 
\begin{enumerate}
\item[(I)] $Q_1$ is finite and $\phi(Q_1) =1$. Equation $\phi'(t) =0$ has only one solution $T_0 \in  [Q_0, Q_1]$ which is the absolute maximum point of $\phi$ on $[Q_0, Q_1]$. Next, if 
 $a> T_0$ is the leftmost point where $\phi'(a) =0$ then $a > Q_1, a -h \in (Q_0, Q_1)$. 
\item[(II)] $\phi$ strongly increases on $(Q_0, +\infty)$, with at most one critical point $Q_0+h$. 
\item[(III)] $\phi$ has exactly two critical points: strong local maximum at \mbox{$T_0 \in (Q_0, Q_0+h)$} and a strong 
local minimum at  $t_m  > Q_0+h$, where $\phi(t_m) \geq 1$. On the interval $(t_m,+\infty)$, function $\phi$ is  increasing with $\phi'(t) >0, \phi''(t)>0$.  
\end{enumerate}
\end{lem}
\begin{pf} \hspace{-2mm} Obviously, we get the second option if $\phi'(t) >0$ for all $t \in \R$. Thus we may suppose that there exists  some leftmost point  $T_0>Q_0$ where $\phi'(T_0)=0$. This implies immediately that    $\phi(T_0) >1$, $\phi''(T_0)\leq 0,$
and, consequently, $\phi(T_0-h) \leq 1$. 

\vspace{-1mm}

(I) Suppose that $Q_1$ is finite so that $\phi(Q_1)=1$. 
We claim that $\phi(T_0-h) <1$ and therefore $\phi''(T_0)<0$ with $T_0$ being a local maximum point.  
Indeed, if $\phi(T_0-h) =1, \phi''(T_0)=0, \phi'(T_0)=0,$ then  $\phi'''(T_0) = \phi(T_0)\phi'(T_0-h) > 0$ in virtue of Lemma \ref{o}. 
In consequence $\phi(t)= \phi(T_0) + \phi'''(T_0)(t-T_0)^3/6 + o((t-T_0)^3)$ and thus $T_0$ is not an absolute maximum 
point on $[Q_0,Q_1]$. Let $q>T_0$ be such a point, then $\phi(q-h) >1, \phi''(q)\leq 0, \phi'(q)=0,$ a contradiction.

\vspace{-1mm}

Hence, $\phi''(T_0)<0, \phi(T_0-h) <1.$ Let $a>T_0$ be the leftmost point where $\phi'(a) =0$. Then $a$ is finite, $\phi''(a) \geq 0$ and therefore $\phi(a-h) \geq 1$.  Now, if $\phi''(a) =0$ then $\phi(a-h) = 1$ and $\phi'''(a) = \phi(a)\phi'(a-h) >0$, a 
contradiction (since $\phi$ is strictly decreasing on $(T_0,a)$).  This means that  $\phi''(a) > 0$ and  $\phi(a-h) > 1$.

\vspace{-1mm}

Suppose that $a < Q_1$.  Then there is $b \in (a, Q_1)$ such that $\phi(b-h)>1,$ $ \phi'(b) =0,$ $ \phi''(b) \leq 0$, 
contradicting to equation (\ref{twe2a}). Next, if $a = Q_1$ then  $\phi(a-h)>1,$ $ \phi(a) =1, $ $ \phi'(a) =0, \phi''(a) > 0$.  Therefore $\phi''(t) >0, t \geq a$, so that  the option (III) holds. (Indeed, otherwise there is $d >a$ such that $\phi'(t) >0$ for all $t \in (a,d],$ and $\phi''(d) =0, \phi(d-h) \geq 1$, a contradiction).

\vspace{-2mm}

(II) Now, suppose that $Q_1 =+\infty$ and $\phi(T_0-h)=1$. Then, as it was shown in (I), we obtain $\phi'''(T_0)>0$ that yields $\phi'(t)>0$ for all $t >T_0$.

\vspace{-2mm}

(III) Finally, consider the situation when   $Q_1 =+\infty$ and $\phi(T_0-h)< 1$ (i.e. $\phi$ reaches a strict local maximum at $T_0$). In such a case, $\phi$ should have subsequent leftmost critical point $q>T_0, \phi(q)>1.$ Indeed, otherwise $\phi'(t)<0, \phi(t) >1,$ $ t > T_0,$ so that $\phi$ converges monotonically to $1$ at $+\infty$.  However, due to the proof of \cite[Lemma 20]{GT}, this is possible only when $(\tau,c) \in \mathfrak{D}_{m}$ and therefore this contradicts to Theorem \ref{ama}.  By the arguments in (I), we already know that $\phi'(q)=0, \phi''(q) > 0$ and  $\phi(q-h) > 1$.  This makes impossible the existence of $p>q$, where  $\phi''(p)=0, \phi'(p) > 0$ and $\phi(p-h)>1$.  In particular, $\phi'(t) >0, t > q.$ \hfill $\square$
\end{pf}

\vspace{-7mm}

\begin{cor} \label{c17} Let $Q_0 < T_0$ be as in  Lemma \ref{mp}(I) or \ref{mp}(III). Then $$\phi(T_0) = \max_{s \in [Q_0,Q_0+h]} \phi(s) \leq e^{ch}$$ and 
$\phi(t) > e^{c(t-Q_0)}, t < Q_0, \phi'(Q_0) < c,   \phi(t) < e^{c(t-Q_0)}, t \in (Q_0,Q_0+h]$. 
\end{cor}

\vspace{-5mm}

\begin{pf}
Integrating equation (\ref{twe2a}) between $-\infty$ and $t \leq Q_0+h$, and taking into account that $\phi(t)(1-\phi(t-h)) >0$ for all 
$t < Q_0+h$, we find that $\phi'(t) - c\phi(t) <0,$ $ t < Q_0+h$. Hence $(\phi(t)e^{-ct})'$ is strictly decreasing on $(-\infty,Q_0+h]$. 
In particular, $\phi'(Q_0) < c \phi(Q_0) =c$ and $\phi(T_0)e^{-cT_0} < \phi(Q_0)e^{-cQ_0} = e^{-cQ_0} $. Thus 
$\phi(T_0) < e^{c(T_0-Q_0)} < e^{ch}.$ The proof of other inequalities is similar. 
\hfill $\square$ \end{pf}
\begin{lem} \label{mpq} Assume that option (I) of  Lemma \ref{mp} holds.  Then there exists a finite number $Q_2 >Q_1$ such that $\phi(Q_2)=1, \phi'(Q_2)>0$ and $\phi(t) <1$ on $(Q_1,Q_2)$. Moreover,  equation $\phi'(t) =0$ has only one solution $T_1 \in  [Q_1, Q_2]$ which is the absolute minimum point on $(Q_1, Q_2)$. Next, if 
 $T_2>T_1$ is the finite leftmost point where $\phi'(T_2) =0$ then $T_2 -h \in (Q_1, Q_2)$. Finally, $Q_2-Q_0 >h$. 
\end{lem}

\vspace{-6mm}

\begin{pf}  Let $T_1>Q_1$ be the leftmost point  where $\phi'(T_1)=0$. By Lemma \ref{mp}(I), we know that  $\phi''(T_1) >0, \ \phi(T_1-h)> 1, \ \phi(T_1) <1$.  Next,  let $(Q_1,Q_2)$ denote  the maximal open interval containing $T_1$  where $\phi(t)<1$.  

\vspace{-2mm}

First, assume that $\phi'(t)>0$ for $t>T_1$. Then $\phi(t)$ is unbounded since otherwise $\phi(t)$ converges monotonically to 1 that is possible only when $(\tau,c) \in \mathfrak{D}_{m}$ and therefore this contradicts to Theorem \ref{ama}. As a consequence, there exists a finite $Q_2$ with the mentioned properties. 

\vspace{-2mm}

Suppose now that there exists some leftmost point $T_2>T_1$ where $\phi'(T_2) =0$.  Then $\phi''(T_2) \leq 0$ and therefore $\phi(T_2-h) \leq 1$.  For an instance, suppose additionally  that  $T_2\in (T_1,Q_2]$.
If $\phi''(T_2) =0$ then $\phi(T_2-h) = 1$ and $\phi'''(T_2) = \phi(T_2)\phi'(T_2-h) < 0$, a 
contradiction (since $\phi$ is strictly increasing on $(T_1,T_2)$).  This means that  $\phi''(T_2) < 0$ and  $\phi(T_2-h) < 1$.
But then  $\phi$ can not have any critical point $b>T_2, \phi(b)<1,$ since otherwise we get  a contradiction: $\phi''(d)=0, \phi'(d) < 0, $ $\phi(d-h) \leq 1$ for some $d \in (T_2,b)$. Therefore $\phi'(t) <0$ for $t > T_2$ so that  
 $\phi''(t) <0$ for $t >T_2$ and $\phi(t)$ can not be positive for large positive $t$. The latter contradiction shows that actually $T_2>Q_2$ and thus $Q_2$ is finite and $\phi'(Q_2)>0$.  Finally, $Q_2-Q_0 > T_1- Q_0 >h$ while the inequality $T_2-h< Q_2$ can be proved in the same way as the inequality $T_0-h <Q_0$ in Lemma \ref{mp}(I).  \hfill $\square$
\end{pf}

\vspace{-6mm}

\begin{cor} Graph of each oscillating solution consists from the arcs similar to described in Lemmas \ref{mp}(I),\ref{mpq} and therefore  it is sine-like slowly oscillating.  
\end{cor}

\vspace{-2mm}

Finally,  the following result describes behavior of positive unbounded waves: 
\begin{cor} \label{c17a} Let profile $\phi$ be unbounded.  Then, for some $A>0$ and $ t_0 \geq Q_0,$ it holds that  $\phi'(t) >0, $ $ \phi(t) > Ae^{ct},$ $ t \geq t_0$.
\end{cor}  

\vspace{-6mm}

\begin{pf} By Lemmas \ref{mp}, parts (II) and (III), for an appropriate $t_0$, each unbounded solution satisfies $\phi'(t) >0, \phi(t-h) >1, \ t \geq t_0$.  This implies that $\phi''(t) - c\phi'(t) >0,\ t \geq t_0$ and therefore 
$\phi'(t) > \phi'(t_0)e^{c(t-t_0)} >0, \ t \geq t_0$. 
\hfill $\square$
\end{pf}

\vspace{-10mm}

\section{{\it A priori} estimates and  the convergence of  semi-wavefronts} \label{ETW}

\vspace{-5mm}

With the change of variables $\phi(t) = e^{-x(t)}$,  equation (\ref{twe2a}) is transformed into  
\begin{equation} \label{et}
x''(t) - cx'(t) - (x'(t))^2 +(e^{-x(t-h)}-1)=0, \ t \in \R. 
\end{equation}

\vspace{-6mm}

Let $\phi(t) = e^{-x(t)}$ be an oscillating semi-wavefront  and for the simplicity take $Q_0=0$. By Corollary {\ref{c17}},   
$
0< x(t) < -ct,  \ t  <0,    
$
and $x(t) > -ct> -ch$ for $t \in (0,T_0)$.  

Our  {\it a priori} estimates are based on the following key assertion: 

\vspace{0mm}

\begin{lem}\label{L20}
Let $y$ solve the boundary value problem 
$$
y'-cy -y^2 + g(t) =0,\  y(a)=y(b)=0,  \  -1 < A:= \min_{s\in [a,b]} g(s) <0,
$$
where $c \geq 2$ and $g$ is continuous. 
Then 
$$
\beta:= \min_{s\in [a,b]} y(s) \geq \frac{2A}{c+\sqrt{c^2+4A}}=:f(A). 
$$
Similarly, 
$$
\gamma:= \max_{s\in [a,b]} y(s) \leq \frac{2B}{c+\sqrt{c^2+4B}} =f(B), \ {\rm where \  }B:=  \max_{s\in [a,b]} g(s). 
$$
\end{lem}

\vspace{-10mm}

\begin{pf}
If $\beta = 0$, the conclusion of the first part of Lemma \ref{L20} is obvious.  Thus we can suppose that $\beta=:y(s')  <0$ and that  $y(t) <0$ for all $t$ from some maximal open interval $(a',b') \subset (a,b)$ containing $s'$.  In particular,  $ y'(s') =0,$ $y(a')=y(b')=0,$ so that 
$\beta \in \{\lambda_1(s'), \lambda_2(s')\}$  where $\lambda_1(s)<  f(A) \leq \lambda_2(s)$ are simple roots of the quadratic equation $y^2+cy - g(s) =0$. Observe that
$$
f(A) \leq \lambda_2(s)= f(g(s)) \leq f(B), 
$$
since $f(u)= 0.5(\sqrt{c^2+4u}-c)$ is strictly increasing in $u$. 
It is clear that each $\lambda_j(s)$ depends continuously on  $s$, and  that 
$\lambda_1(s) < f(A) <0$. Suppose for a moment that $\beta = \lambda_1(s')$. We claim that then $y(s) < f(A)$ for all $s \in [s',b']$. Indeed, let $q$ be the minimal real number such that  $y(q) = f(A)$. 
Then $y'(q) \geq 0$ and we have the following dichotomy: either (i) $y'(t) >0$ on some maximal subinterval $(p,q),\  y'(p) =0$, of $(s',q)$, or (ii) 
there exists a sequence $\{t_j\},\  t_j <q,$ converging to $q$ such that $y'(t_j)=0$. In every case,  $y(p) = \lambda_1(p)$,  $y(t_j) = \lambda_1(t_j)$ due to $y(p), y(t_j) < y(q) = f(A)$.  Therefore  the case (i) is not possible because of the following contradiction: 
$y(q)= \lim y(t_j) = \lim \lambda_1(t_j) = \lambda_1(q) < f(A)$.  Similarly, the case (ii) should also be discarded
in virtue of the following argument:  $y^2(s)+cy(s)-g(s) = y'(s) >0$ on $(p,q)$, $y(p) = \lambda_1(p),$ so that $y(s) < \lambda_1(s) < f(A), \ s \in (p,q)$, whence $y(q) \leq  \lambda_1(q) < f(A)$,  a contradiction.    
Hence, assuming that 
$\beta = \lambda_1(s')$  we obtain that  $y(s) < f(A)$ for all $s \in [s',b']$. In particular, 
 $0= y(b') < f(A) <0$. This contradiction  proves  the first part of Lemma \ref{L20}. 
 
Next, it is clear that $\gamma \geq 0$.  If $\gamma =0$ then $B \geq g(a) = - y'(a) \geq 0$ and the claimed inequality is immediate. If $\gamma >0$ then $\gamma \in \{\lambda_1(s''), \lambda_2(s'')\}$  for some $s'' \in (a,b)$. As a consequence, we obtain the second estimation of the lemma:  
$
\lambda_1(s'') \leq \gamma \leq  \lambda_2(s'') \leq f(B). 
$
\hfill $\square$
\end{pf}   

Recall  that the  Schwarz derivative $Sp$ of $C^3$-smooth function $p$ is defined as
$$(Sp)(x)=p'''(x)(p'(x))^{-1}-(3/2) \left(p''(x)(p'(x))^{-1}\right)^2.$$
\begin{lem} \label{nsw} Let $c \geq 2$. Then  real analytic  function 
$(f\circ g)(x)= f(e^{-x}-1), $ $ x \in \R, $ $ (f\circ g)(0) =0,$ is well defined, strictly decreasing and has the 
negative Schwarz derivative on $\R$. 
\end{lem}

\vspace{-6mm}

\begin{pf}   Since $f(u)= 0.5(\sqrt{1+4u/c^2} -1)$, 
we find easily that $(Sf)(u)= 6(c^2+4u)^{-2}$.
By  the well known formula 
for the Schwarzian of the composition, 

\vspace{-10mm}

$$
S(f\circ g)(x) = (Sf)(g(x))(g'(x))^2 + (Sg)(x) = \frac{6e^{-2x}}{(c^2+4(e^{-x}-1))^2}- \frac12= 
$$
$$
= \frac{6}{(e^x(c^2-4)+4)^2}- \frac12 \leq \frac{6}{4^2} -\frac12 =-\frac18.
$$
The other properties of $f\circ g$ are straightforward to verify. \hfill $\square$
\end{pf}
\vspace{-5mm}
\begin{lem} \label{ogran} Let $c\geq 2$ and $\phi(t), \phi(-\infty) =0,$ be a slowly oscillating on $[Q_0, +\infty)$ positive solution of  equation 
(\ref{twe2a}). Then $\phi$ is bounded and
$$
0<L_e(c,h) < \phi(t) < U_e(c,h), \ t \geq Q_0, 
$$
where $U_e(c,h):= \exp(-L(c,h)), \ L_e(c,h):= \exp(-U(c,h))$ and  $$U(c,h)= hf(e^{-L(c,h)}-1), \quad 
L(c,h):= -h\max\left\{c, \frac{2}{c+\sqrt{c^2-4}}\right\}.
$$
\end{lem}

\vspace{-10mm}

\begin{pf} Without the loss of generality, we can set $Q_0=0$. Then it suffices to prove the boundedness of $x(t) = -\ln \phi(t)$ on $[0,+\infty)$. Since $\phi(t)$ is  slowly oscillating about $1$, the transformed solution  $x(t)$ oscillates slowly around the zero equilibrium of (\ref{et}).  This implies that 
there exists an increasing  sequence $Q_j, j \geq 0, Q_0=0,$ of zeros of $x(t)$ such that $x(t)<0$ on 
$(Q_0,Q_1)\cup(Q_2,Q_3) \cup \dots$ and $x(t)>0$ on 
$(Q_1,Q_2)\cup(Q_3,Q_4) \cup \dots$ We proceed by evaluating extremal values $V_j=x(T_j)$ of $x(t)$ on each  interval $(Q_j,Q_{j+1})$.  As we already have established, $|V_0| = -V_0 \leq ch$. 
Next,   we have that $V_1= x(T_1) >0$ with $T_1 > h$ and $x'(T_1) =0, \  x(Q_1) =0, \ T_1-Q_1 <h$. Hence, 
$$V_1 =  \int_{Q_1}^{T_1}x'(s)ds\leq h \max_{s \in [T_0,T_1]} x'(s) \leq h\max_{s \in [T_0,T_1]}f(e^{-x(s-h)}-1) \leq hf(w(V_0)),$$
where $w(x):= e^{-x}-1$.
Next, consider $V_2= x(T_2) <0$, we have   $x'(T_2) =0,$ $ x(Q_2) =0$ and $T_2- Q_2<h$.  Recalling that $\phi(t)$ (and, consequently, $x(t)$) is sine-like slowly oscillating (so that $x'(t) <0$ on $(T_1,T_2)$) and  applying Lemma  \ref{L20}, we obtain
$$
V_2 =  \int_{Q_2}^{T_2}x'(s)ds\geq  h \min_{s \in [T_1,T_2]} x'(s) \geq 
h\min_{s \in [T_1,T_2]}f(e^{-x(s-h)}-1) \geq hf(w(V_1)).$$ 
Since $Q_{j+2}-Q_j >h$ for each $j$, we can repeat the above two steps  to conclude that 
\begin{equation}\label{raz}
V_{2j+1} \leq hf(w(V_{2j})), \ j \geq 0, \quad V_{2j} \geq hf(w(V_{2j-1})), \ j >0. 
\end{equation}

\vspace{-3mm}

As a consequence,

\vspace{-3mm}

$$
V_{2j} \geq hf(w(V_{2j-1})) > hf(w(+\infty)) =   \frac{-2h}{c+\sqrt{c^2-4}}=:B_*(c,h), \ j >0, 
$$
and therefore, after setting $L(c,h)= \min\left\{-ch,B_*(c,h)\right\}$, we obtain that 
$$
V_{2j+1} \leq hf(w(V_{2j})) \leq hf(w(L(c,h))), \ j \geq 0. 
$$
This ends the proof of Lemma \ref{ogran}. \hfill $\square$
\end{pf}

\vspace{-6mm}

Suppose now that $\phi, \ \phi(-\infty)=0,$ is an unbounded positive solution of (\ref{twe2a}). By Lemmas \ref{mp}, \ref{mpq} and \ref{ogran}, function $\phi$ is either monotone or slowly oscillating around $1$ on some interval $(-\infty, Q_m], \ \phi(Q_m) =1,$ and $\phi(t) >1$ for $t > Q_m$. Let $T_m$ denote the rightmost critical point of 
$\phi$ (whenever it exists) and set $S_m = \max\{T_m,Q_m\}$. 

\begin{cor} \label{c22or}There exists a positive constant $\beta(c,h)> U_e(c,h)$ depending only on $(c,h)$ such that 
$\phi(t) < \beta(c,h), \ t \leq S_m+2h$.  In this way, if $\phi(\bar s) = \beta(c,h)$ for some $\bar s\in \R$ then  $\phi'(t) >0,\ \phi(t) >1$ for all $t \geq \bar s-h$.
\end{cor}

\vspace{-6mm}

\begin{pf} \underline{Step I.}  Suppose first that $S_m= Q_m> T_m$.
As the proof of Lemma \ref{ogran} shows, we have that $\phi(t) < U_e(c,h)$ for all $t \leq Q_m$.  Next, 
on the half-line  $\mathcal{I}:=(-\infty,Q_m+h]$, function $\phi(t)$ satisfies the homogeneous linear equation 
\begin{equation}\label{gle}
y''(t) - cy'(t) + a(t)y(t) =0, \end{equation}

\vspace{-3mm}

whose coefficient $a(t) = 1-\phi(t-h)$  is uniformly bounded on $\mathcal{I}$ by a constant depending only on  $c,h$.  We consider separately the cases $m=0$ and $m>0$.  

If $m=0$ then $\phi'(t) >0$ for all $t$, and 
$\phi'(Q_0) < c, \phi(Q_0)=1$, see Corollary \ref{c17}.  As a consequence, the solution $y(t)\equiv \phi(t)$ of the initial value 
problem $y(Q_0)= 1, y'(Q_0) = \phi'(Q_0),$ to equation  (\ref{gle}) exists on $(-\infty,Q_0+h]$ where it is  bounded  by some constant $\rho_0(c,h)$ depending only on $c,h$.  Therefore the absolute value of $a(t) =1- \phi(t-h) = 1-y(t-h), \ t \leq Q_0+2h,$ is bounded by $\rho_0(c,h)+1$ and we can repeat the above argument to conclude that  the solution $y(t)\equiv \phi(t)$ of the initial value 
problem $y(Q_0)= 1, y'(Q_0) = \phi'(Q_0),$ of equation  (\ref{gle}) exists on $(-\infty,Q_0+2h]$ where it is  bounded by some constant $\rho_1(c,h)$.

Now we can assume that $m>0$ and $\phi(t) <1$  on some maximal interval $(Q_{m-1},Q_m)$.  We also know that $\phi'(t) >0$ on some maximal open subinterval $(T_{m-1},Q_m)$
of $(Q_{m-1},Q_m)$.  Since $\phi'(T_{m-1})=0, \phi(T_{m-1}) <1$, we can  integrate equation (\ref{gle})  repeatedly (as it has been done in the case $m=0$)   to prove the existence of  $\rho_2= \rho_2(c,h)$ such that 
$
\phi(t) < \rho_2, \ t \leq T_{m-1}+4h. 
$
If $Q_m+2h \leq T_{m-1}+4h$, the proof is finished. Otherwise $Q_{m} > T_{m-1}+2h$ and $\phi(t)$ is strictly increasing 
on $[Q_m-2h,Q_m]$. In particular, $\phi'(\hat s) \in (0,  (2h)^{-1}),$ $ \phi(\hat s) \in (0,  1),$ at some point $\hat s \in [Q_m-2h,Q_m]$. But then there exists $\rho_3 = \rho_3(c,h)$ depending only on $c,h$ and such that 
$\phi(t) < \rho_3$ on $[\hat s,\hat s+4h] \supset [Q_m,Q_m+2h]$.  Therefore, by taking $\beta(c,h) = \max\{\rho_j(c,h), j =0,1,2,3\}$, we finalize  the proof of  Corollary \ref{c22or} in the case when $S_m= Q_m> T_m$.

\underline{Step II.}  Suppose now that $S_m= T_m\geq Q_m$. This situation corresponds to the cases $(II)$ and $(III)$ of Lemma \ref{mp}.  Since $\phi'(S_m) =0$ and $\phi(t) \leq U_e(c,h),$ $ t \leq S_m$,  we can again  integrate equation (\ref{gle})  repeatedly  to prove the existence of  $\rho_4= \rho_4(c,h)$ such that 
$
\phi(t) < \rho_4, \ t \leq S_{m}+2h. 
$ \hfill $\square$
\end{pf}

\vspace{-6mm}

For fixed $c\geq 2,h>0$, we will consider also the following modified equation 
\begin{equation}
\label{twe2m} 
\phi''(t) - c\phi'(t) + g(\phi(t))(1-
\phi(t-h)) =0,  
\end{equation}
with $\beta(c,h)$ defined  in Corollary \ref{c22or} and with continuous piece-wise linear
$$g(u)=\left\{\begin{array}{cc} u,& u \in [0,\beta(c,h)], \\  \max\{0,2\beta(c,h) -u\},
& u>  \beta(c,h).\end{array}\right. $$ 
\begin{lem}
\label{twoe} Equations (\ref{twe2m}) and (\ref{twe2a}) share the same set of semi-wavefronts. 
\end{lem}

\vspace{-6mm}

\begin{pf} Due to Lemma \ref{ogran} and the definition of $g(u)$, each  semi-wavefront of 
(\ref{twe2a})  also satisfies (\ref{twe2m}).   Conversely, suppose that $\phi$ is a semi-wavefront to 
(\ref{twe2m}).  We will prove that then $\phi(t) < \beta(c,h)$.  Indeed, otherwise  $\phi(\bar s) = \beta(c,h)$ at some leftmost  point $\bar s$.  Since $\phi(t)$ is also satisfying  (\ref{twe2a}) for all $t\leq \bar s$, Corollary \ref{c22or} assures that 
$\phi'(\bar s) >0$ and $\phi(t-h) > 1$ for all $t \in [\bar s,\bar s+h]$. Thus $\phi''(t) >  c\phi'(t),$ $t \in [\bar s,\bar s+h],$ and  consequently 
$\phi'(t)> \phi'(\bar s)e^{c(t-\bar s)}, t \in [\bar s,\bar s+h]$. Hence, 
$\phi''(t) \geq c\phi'(t) >0$ on $[\bar s,\bar s+h]$. Using step by step continuation argument, we can conclude that $\phi(+\infty)=+\infty$, a contradiction. \hfill $\square$
\end{pf}
\begin{lem}
\label{ocov} Let $\phi(t)$ be a slowly oscillating semi-wavefront to equation 
(\ref{twe2a}). If $\tau \leq 1, \ c \geq 2$, then $\phi(+\infty)=1$.
\end{lem} 
\begin{pf} Lemma \ref{ogran} assures the existence of  finite limits 
$$0 \geq m_*= \liminf_{j \to +\infty} V_j = \liminf_{t \to +\infty}x(t), \quad 0 \leq  M_* = \limsup_{j \to +\infty} V_j=
\limsup_{t \to +\infty}x(t).$$ 
Clearly, the lemma will be proved if we show that $\tau \leq 1$ implies $M_*=0$.
From (\ref{raz}), we deduce that 
$
M_* \leq hf(w(m_*)), \  m_* \geq hf(w(M_*))
$ and therefore 
$$
M_*\leq (hf\circ w)^2(M_*) \leq  (hf\circ w)^4(M_*)\leq \dots \leq (hf\circ w)^{2k}(M_*) \leq \dots 
$$
Here $f^k = f\circ \dots \circ f$ denotes the $k$-times  composition of $f$.  Now, by Lemma \ref{nsw}, analytic function  $h f\circ w$ is strictly decreasing, below bounded and has the negative Schwarzian.  Therefore the inequality $|hf'(0)w'(0)| = h/c = \tau \leq 1$ assures the global stability of the fixed point $0$ of the one-dimensional mapping $h f\circ w: \R \to \R$.  See \cite[Proposition 3.3]{ltt} for more details.  In particular, this means that $(hf\circ w)^{2k}(M_*) \to 0$ as $k \to +\infty$.  Hence, $M_*=0$ and Lemma \ref{ocov} is proved. 
\hfill $\square$
\end{pf}

\section{Existence of semi-wavefronts for $c\geq 2, \ h >0$.} \label{ESW}
With $g(u)$ defined in (\ref{twe2m}) and with some $b>1+2\beta(c,h)$ (to be specified later), let us consider   $r(\phi(t),\phi(t-h)):= b\phi(t) + g(\phi(t))(1-
\phi(t-h))$.  By Lemma \ref{twoe}, it suffices to prove that equation  
\begin{equation}
\label{twe2mm} 
\phi''(t) - c\phi'(t)  - b \phi(t) + r(\phi(t),\phi(t-h))=0
\end{equation}
has a semi-wavefront.
Observe that if some  $\psi(t)$ satisfies $0\leq \psi(t) \leq \beta(c,h)$ and $\psi(t-h) \leq 2\beta(c,h) <b$, then  
\begin{equation}\label{Ar}
r(\psi(t),\psi(t-h))=  \psi(t)(b +1-
\psi(t-h)) \geq 0. 
\end{equation}
 Now, if   $\beta(c,h)\leq \psi(t) \leq 2\beta(c,h)$ and $\psi(t-h) \leq 2\beta(c,h) <b$, then  $$
r(\psi(t),\psi(t-h))=  (2\beta(c,h) -\psi(t))(1-\psi(t-h))+ b
\psi(t) =
$$
\begin{equation} \label{Br}
2\beta(c,h)(1-\psi(t-h))+ \psi(t) (b-1+\psi(t-h))> \beta(c,h).  
\end{equation}
Next, we consider the non-delayed KPP-Fisher equation  $u_t= u_{xx} +g(u)$.  The profiles $\phi$
of the traveling fronts $u(x,t)= \phi (x+ct)$ for this equation satisfy
\begin{equation}\label{kppm}
\phi''(t) - c\phi'(t)  + g(\phi(t))=0, \ c \geq 2.
\end{equation}
As before,  $0< \lambda \leq \mu$ denote eigenvalues of equation (\ref{kppm}) linearized around $0$.  Then 
$\chi(\lambda) = \chi(\mu) =0$ where $\chi(z) := z^2-cz+1$.  Recall also that $z_1< 0< z_2$ stand for the  roots of the equation  $z^2 -cz-b =0$. 
In the sequel,  $\phi_+(t)$ will denote the unique monotone front to (\ref{kppm}) normalized by the condition 
$$ \phi_+(t):= (-t)^je^{\lambda t} (1+o(1)), \ t \to -\infty. $$ In fact,  the latter asymptotic formula can be considerably improved since $\phi_+(t)$ for all $t$ such that $\phi_+(t) < \beta(c,h)$ satisfies the linear differential equation 
$$
\phi''(t) - c\phi'(t)  + \phi(t)=0.  
$$
For example, if $c >2$  then there exists (cf. \cite[Theorem 6]{GT}) $K \geq 0$ such that  
\begin{equation} \label{1mo}
 \phi_+(t):= e^{\lambda t}  - K e^{\mu t}, \ t \leq \phi_+^{-1}(\beta(c,h)).   \end{equation}
 \vspace{-3mm}
Set $\epsilon' =z_2-z_1$  and consider the following integral operator 
$$(A_m\phi)(t) = \frac{1}{\epsilon'}\left\{\int_{-\infty}^te^{z_1
(t-s)}r(\phi(s), \phi(s-h))ds + \int_t^{+\infty}e^{z_2
(t-s)}r(\phi(s), \phi(s-h))ds \right\}.
$$
\begin{lem}  \label{usl} Assume that  $b > 2\beta(c,h) +1$ and let  $0\leq \phi(t) \leq \phi_+(t)$, then $\phi_+$ is an upper solution:  
$$0 \leq (A_m\phi)(t) \leq \phi_+(t).$$
\end{lem}

\vspace{-5mm}

\begin{pf} The lower estimate is obvious since $0\leq \phi(t) \leq \phi_+(t) \leq 2\beta(c,h)$ and therefore
$r(\phi(t),\phi(t-h)) \geq 0$ in view of (\ref{Ar}) and  (\ref{Br}). Now, since  $\phi(t) \leq \phi_+(t)$ and $bu+g(u)$ is an increasing function, we find that 
$$
r(\phi(s), \phi(s-h)) \leq b\phi(t) + g(\phi(t)) \leq b\phi_+(t) + g(\phi_+(t)) =:R(\phi_+(t)). 
$$
Thus 
$$
(A_m\phi)(t) \leq  \frac{1}{\epsilon'}\left\{\int_{-\infty}^te^{z_1
(t-s)}R(\phi_+(s))ds + \int_t^{+\infty}e^{z_2
(t-s)}R(\phi_+(s))ds \right\}= \phi_+(t),
$$
and the lemma is proved. \hfill  $\square$ 
\end{pf}

\vspace{-5mm}

Next,  we need to find a lower solution for (\ref{twe2mm}).  Fortunately, for $c >2$ we can use the following  well known  solution (e.g. see \cite{wz})  
$$
\phi_-(t)= \max\{0,e^{\lambda t} (1- Me^{\epsilon t})\},
$$
where $\epsilon \in (0, \lambda)$  and $M \gg 1$ is chosen in 
such a way that $- \chi(\lambda+\epsilon)   > 1/M,$ $\lambda +\epsilon < \mu$, and
$$0< \phi_-(t) < \phi_+(t) < e^{\epsilon t}< 1, \quad t \leq T_c, \ \mbox{where} \ \phi_-(T_c) =0.$$
The above inequality $\phi_-(t) < \phi_+(t)$ is possible due to  representation (\ref{1mo}). 
\begin{lem}  \label{nre} Let $b > 2\beta(c,h) +2$ and $\phi_-(t) \leq \phi(t) \leq \phi_+(t), \ t \in \R$, then 
\begin{equation}\label{ul2}
\phi_-(t) \leq (A_m\phi)(t) \leq \phi_+(t), \quad t \in \R. 
\end{equation}
\end{lem}

\vspace{-8mm}

\begin{pf} Due to Lemma \ref{usl}, it suffices to prove the first inequality in (\ref{ul2}) for $t \leq T_c$. 
Since $0< \phi(t) < 1 < \beta(c,h), \ t \leq T_c$, we have, for $t\leq T_c $ that 
$$
(A_m\phi)(t) \geq   \frac{1}{\epsilon'}\left\{\int_{-\infty}^te^{z_1
(t-s)}r(\phi(s), \phi(s-h))ds + \int_t^{T_c}e^{z_2
(t-s)}r(\phi(s), \phi(s-h))ds \right\}=
$$
$$
\frac{1}{\epsilon'}\left\{\int_{-\infty}^te^{z_1
(t-s)} \phi(s)(b +1-
\phi(s-h)) ds + \int_t^{T_c}e^{z_2
(t-s)}\phi(s)(b +1-
\phi(s-h)) ds \right\}\geq 
$$
$$ \frac{1}{\epsilon'}\left\{\int_{-\infty}^te^{z_1
(t-s)} \phi_-(s)(b +1-
\phi_+(s-h)) ds + \int_t^{T_c}e^{z_2
(t-s)}(\dots) ds \right\}= 
$$
$$ \frac{1}{\epsilon'}\left\{\int_{-\infty}^te^{z_1
(t-s)} \phi_-(s)(b +1-
\phi_+(s-h)) ds + \int_t^{+\infty}e^{z_2
(t-s)}(\dots) ds \right\} =: Q(t),
$$
where $(\dots)$ stands for $\phi_-(s)(b +1-
\phi_+(s-h))$. 
In order to evaluate $Q(t)$, we consider
the following chain of inequalities (for $t \leq T_c$)
$$
\phi_-''(t) - c\phi_-'(t)  -b\phi_-(t) + b\phi_-(t) + \phi_-(t)(1- \phi_+(t-h)) =
$$
$$
- \chi(\lambda+\epsilon) M e^{(\lambda +\epsilon)t} - \phi_+(t-h)e^{\lambda t} (1- Me^{\epsilon t})\geq 
$$
$$
- \chi(\lambda+\epsilon) M e^{(\lambda +\epsilon)t} - e^{\epsilon(t-h)}e^{\lambda t} > M e^{(\lambda +\epsilon)t}(
- \chi(\lambda+\epsilon)   - 1/M) >0. 
$$
But then, rewriting the latter differential inequality in the equivalent  integral form  (e.g. see \cite{ma}) and using  the fact that 
$$\Delta\phi'_-|_{T_c} := \phi'_-(T_c+) - \phi'_-(T_c-) = - \phi'_-(T_c-) >0, $$
we may conclude that $Q(t) \geq \phi_-(t), \ t \in \R$.  Hence, $(A_m\phi)(t) \geq  \phi_-(t), \ t \in \R$, and Lemma \ref{nre} is proved. \hfill $\square$
\end{pf}

\vspace{-5mm}

Finally, it is clear that,  in order to establish the existence of semi-wavefronts to equation (\ref{twe2mm}), it suffices to prove that 
the equation $A_m\phi=\phi$ has at least one  solution  from the subset 
$K = \{x \in X: \phi_-(t) \leq x(t) \leq  \phi_+(t), \ t \in \R\}$
of the Banach space $(X,\|\cdot\|)$, where 

\vspace{-5mm}

\begin{eqnarray*}
  X &=& \{x \in C(\R,\R): \|x\| =
\sup_{s \leq 0} e^{-\lambda s/2 }|x(s)|+ \sup_{s \geq 0} e^{-\rho
s}|x(s)|<
\infty\}\end{eqnarray*}
is defined  with some fixed $\rho >0$.  Observe that 
the convergence $x_n \to x$ on $K$ is equivalent to the
uniform convergence $x_n \Rightarrow x$ on compact subsets of
$\R$.
\begin{lem} Take $c >2$. Then $K$ is a closed, bounded, convex subset of $X$
and $A_m:K \to K$ is completely continuous.
\end{lem}

\vspace{-6mm}

\begin{pf} By the previous lemma, $A_m(K) \subset K$. It is also obvious that $K$ is a  closed, bounded, convex subset of $X$. Since 
\begin{equation}\label{AA}
|x(t)|+ |(A_mx)'(t)|\leq 2\beta(c,h)(1+\epsilon'),\  \mbox{for all} \   x \in K,  
\end{equation}

due to the Ascoli-Arzel${\rm
\grave{a}}$ theorem  $A_m(K)$ is precompact in $K$ . 
Next, by the Lebesgue's dominated convergence theorem, if $x_j\to  
x_0$ in $K$ then 
$(A_mx_j )(t) \to  
(A_mx_0 )(t)$ at every $t \in  
\R$. The precompactness of 
$\{A_mx_j \} \subset K$ assures that, in 
fact, $A_mx_j\to A_mx_0$ in $K$. Hence, the map $A_m: K\to K$ is completely  continuous. 
\end{pf}

\vspace{-6mm}

\begin{thm} \label{34} Assume that $c\geq 2$. Then the integral equation
$A_m\phi=\phi$  has at least one  positive bounded solution in $K$.
\end{thm}
\begin{pf}
If $c>2$ then, due to the previous  lemma, we can apply the Schauder's fixed point
theorem to $A_m:K \to K$. Let now $c=2$ and consider $c_j:= 2+1/j$ with $h_0:=2\tau, h_j:= c_j\tau$. By the first part of the 
theorem, we know that for each $c_j$ there exists a semi-wavefront $\phi_j$: we can 
normalize it by the condition $\phi_j(0)= 1/2, \ \phi_j'(s) >0, s \leq 0$. It is clear from (\ref{AA}) that the set 
$\{\phi_j, j \geq 0\}$ is precompact in $K$ and therefore we can also assume that $\phi_j \to \phi_0$ in $K$, where 
$\phi_0(0) =1/2$ and $\phi_0$ is monotone increasing on $(-\infty,0]$.  In addition, 
$R_j(s):= r(\phi_j(s), \phi_j(s-h_j)) \to R_0(s):= r(\phi_0(s), \phi_0(s-h_0))$ for each fixed $s \in \R$. The sequence 
$\{R_j(t)\}$ is also uniformly bounded on $\R$. All this allows us to apply the Lebesgue's dominated convergence theorem in  
$$(A_{m,j}\phi_j)(t) := \frac{1}{\epsilon_j'}\left\{\int_{-\infty}^te^{z_{1,j}
(t-s)}R_j(s)ds + \int_t^{+\infty}e^{z_{2,j}
(t-s)}R_j(s)ds \right\} = \phi_j(t), 
$$
where $z_{1,j} <0< z_{2,j}$ satisfy $z^2-c_jz -b =0$. 
In this way we obtain that $A_m\phi_0 = \phi_0$ with $c=2$ and therefore $\phi_0$ is a non-negative solution 
of equation (\ref{twe2a}) satisfying condition $\phi_0(0)=1/2$ and monotone increasing on $(-\infty,0]$. It is immediate to see that $\phi_0(-\infty)=0$ and therefore $\phi_0$ is a semi-wavefront. \hfill $\square$
\end{pf}

\vspace{-10mm}

\section{Admissible wavefront speeds} \label{LaT}

\vspace{-5mm}

First, we observe that the necessity of the condition $c\geq 2$ for the existence of monotone wavefronts was 
already established in \cite[Lemma 19]{GT}. Since the leading edge of each semi-wavefront is  monotone, 
the proof of the mentioned lemma is also valid for the broader family of semi-wavefronts. 

Consider now some semi-wavefront $\phi$ propagating at the velocity $c> c^\star$.   We know that  $\phi$ is slowly oscillating around the positive steady state. In this section, we show that
these oscillations are non-decaying.

Arguing by contradiction,  assume that $\phi(+\infty) =1$.  Then $w(t)=\phi(t)-1,$ $ w(+\infty) =
0,$ solves
$$
w''(t) - cw'(t)-  w(t-h)(1+w(t))=0,\quad  t \in \R.\
$$

\vspace{-5mm}

Since $w(+\infty) = 0$,  there exists a subsequence $\{t_n\},\ \lim t_n=+\infty,$ of the sequence $\{T_n\}$ defined in Lemma \ref{ogran} such that $|w(t_n)| = \max_{s \geq t_n}|w(s)| >0,$ $w'(t_n) =0, $
$w''(t_n)w(t_n) < 0, \ w(t_n)w(t_n-h) <0$.  In fact, there is a unique $q_n \in (t_{n}-h,
t_{n})$  such that $w(q_n)=0$. Without restricting the generality, we can suppose that $w(t_n) >0$ and that 
$\{r_n\}, \
r_n:=t_n- q_n \in (0,h),$ is monotonically converging to $ r_*\in
[0,h]$.  Clearly, $w(s) < 0$ for $s \in
[t_{n}-h, q_{n})$ and $w(s)  > 0$ for $s \in (q_{n},
t_{n}]$.

Now, each $y_n(t) : = w(t+t_n)/w(t_n), \ t \in \R, $ satisfies
\begin{equation}\label{twerin}  \hspace{0mm}
y''(t) - cy'(t)- (1+w(t+t_n))y(t-h)=0. 
\end{equation}
 It is clear that $y_n(0)=1$ and $|y_n(t)| \leq 1, \ t
\geq -r_n.$ In addition,  $ y_n(-r_n)=0,$ $ y_n(-h) <0$. We also can suppose that
$|w(t+t_n)| \leq 0.1$ for all $n$ and $t \geq 0$. 

\vspace{-3mm}

Next, we are going  to estimate $|y'_n(t)|, \ t \geq 0$. Let $s\geq 0$ be the leftmost  local extremum point 
for $y_n'(t)$. Then $y_n''(s)=0, y_n'(s) <0,$ and therefore
$$
0>cy_n'(s) =-  y_n(s-h)(1+w(s+t_n)). 
$$
Thus $y_n(s-h) >0$ that yields $s-h > -r_n$.  Consequently,  $\bar s -h >-r_n$ for each other 
critical point $\bar s$ of $y_n'(t)$. All this implies  that $|y_n(\bar s-h)| \in [0,1]$. Therefore $|y'_n(t)| \leq 1.1/c$ for $t \geq 0$, and, in particular, $y_n(t) \geq 0.45$ on $[0,c/2]$.  
Next, due to the Ascoli-Arzel${\rm
\grave{a}}$ theorem,  the sequence $y_n(t)$ has a
subsequence which converges on $[0, +\infty)$, in the compact-open
topology, to some continuous function $y_*(t)$.  Evidently,  $ \max \{|y_*(s)|, s \geq 0\} = y_*(0)= 
1$ and  $y_*(t) \geq 0.45$ on $[0,c/2]$.   Next, for some fixed positive $b$ and all $t \in [h,+\infty)$, it holds that
$$g_n(t):= by_n(t)- (1+w(t+t_n))y_n(t-h) \to g_*(t):= by_*(t)-y_*(t-h).
$$
Obviously, $0 \leq |g_*(t)| \leq 1+b$ for $t \geq h$.

In order to establish some further properties of
$y_*(t)$, let us present the family of all solutions to $(\ref{twerin})$
which are bounded at $+\infty$:
\begin{equation}\label{rereink} \hspace{0mm}
\hspace{-7mm} y(t) = Ae^{z_1 t} + \frac{1}{\epsilon'}\left\{ \int_{h}^te^{z_1 (t-s)}g_n(s)ds +
\int_t^{+\infty}e^{z_2 (t-s)}g_n(s)ds \right\}, \ t \geq h.
\end{equation}

\vspace{-3mm}

Here $\epsilon' = z_2-z_1$ is defined in the same way as in Lemma \ref{usl}. Replacing $y(t)$ with $y_n(t)$ in (\ref{rereink}), we obtain that, for some $A_n$, 
\begin{displaymath}
\hspace{-7mm} y_n(t) = A_ne^{z_1 t} + \frac{1}{\epsilon'}\left\{ \int_{h}^te^{z_1 (t-s)}g_n(s)ds +
\int_t^{+\infty}e^{z_2 (t-s)}g_n(s)ds \right\}, \ t \geq h.
\end{displaymath}

\vspace{-3mm}

The latter inequality implies that $A_n, \ n \in \N,$ are uniformly  bounded: 
$$
|A_n| = e^{-z_1h} \left | y_n(h)- \frac{1}{\epsilon'}
\int_h^{+\infty}e^{z_2 (h-s)}g_n(s)ds\right| \leq  e^{-z_1h} \left (1+ \frac{1.1+b}{\epsilon'z_2}\right).
$$
Hence, 
taking limit as $n \to + \infty$ (through passing to a subsequence
if necessary) we find that $y_*(t)$ satisfies
\begin{equation}\label{rerein} \hspace{-0mm}
y_*(t) = Ae^{z_1 t} + \frac{1}{\epsilon'}\left\{ \int_{h}^te^{z_1 (t-s)}g_*(s)ds +
\int_t^{+\infty}e^{z_2 (t-s)}g_*(s)ds \right\},\  t \geq h,
\end{equation}
with some finite $A$. Now, (\ref{rerein}) implies that $y_*(t)$
is a solution of  the equation
\begin{equation}\label{eqwo}
y''(t) - cy'(t)- y(t-h)=0, \ t \geq h.
\end{equation}

\vspace{-6mm}

We claim that $y_*(t)$ is not a small solution. Indeed, on the contrary, let us suppose that $y_*(t)$
has superexponential decay. Since the characteristic function 
$z^2-cz-e^{-zh}$ to (\ref{eqwo}) has the exponential type $h$, an application of  
\cite[Theorem 3.1]{hale} assures
that $y_*(t) = 0$ for all $t \geq 2h$. But then equation (\ref{eqwo})
implies that $y_*(t) = 0$ for all $t \geq h$ and, in
consequence, $y_*(t) = 0,$ for all $t \geq 0$. This contradicts the inequality 
$y_*(t) \geq 0.45$ on $[0,c/2]$ and therefore $y_*(t)$ is not a small solution.

Hence, by \cite[Proposition 7.2]{FA}, for every
sufficiently large $\nu < 0$, we have that
$$
y_*(t) =  u(t) + O(\exp(\nu t)), \ t \to +\infty,
$$
where $u$ is a {\it non empty} finite sum of  eigensolutions of
(\ref{eqwo}) associated to the eigenvalues $\lambda_j \in F= \{\nu
< \Re\lambda_j \leq 0\}$. Now, Lemmas \ref{c1de} and \ref{PL} say that, for
every $c > c^\star$,
$$F\cap
(-\infty,0]\times [-2\pi/h, 2\pi/h] = \emptyset.
$$
In consequence, there exist $A >0,\ \beta > 2\pi/h,\ \alpha \geq 0,\
\xi \in \R$, such that
$$
y_*(t) =  (A\cos (\beta t+ \xi) + o(1))e^{-\alpha t}, \ t \geq
0.
$$
This implies the existence of an interval $(a,a+h)$,
$a
> 3h$, such that $y_*(t)$ changes its sign on $(a,a+h)$ exactly three
times. Since $y_{n_j}(t) \to y_*(t)$ uniformly on
$[a,a+h]$, we can conclude that sc$(\bar y_{n_j, a+h})\geq 3$ for
all large $j$. However, this contradicts to  the slowly
oscillating behavior of $y_{n_j}(t)$. In consequence,  the equality  $\phi(+\infty) =1$ can not hold for $c> c^\star$.

\vspace{-3mm}

\section{Appendix: Proof of Lemma \ref{c1de} }

\vspace{-3mm}

In this section, we study the zeros of $\psi(z,c):=
z^2-cz-e^{-z c\tau}, $  $c \geq 2,$ $\tau >0.$ It is straightforward  to see
that $\psi$ always has a unique positive simple zero $\lambda_{-1}$. Since
$\psi'''(z,c)$ is positive, $\psi$ can have at most three
(counting multiplicities) real zeros, one of them positive and the
other two (when they exist) negative.

\vspace{-2mm}

Fix some $\tau \geq 0$. We should prove that $\psi(z,c), c \geq 2,$ has exactly one (counting multiplicity) 
zero in the open right half-plane $\Re z > 0$ if and only if $\tau \leq \tau_2$ and $c \leq c^\star(\tau)$.  Aiming this objective, we first consider $\psi(z,c)$ without restriction $c \geq 2$. Then the next result follows from \cite[Section 2]{GT}:

\vspace{-2mm}

\begin{lem}  \label{roots} There exists function $C^*=C^*(\cdot): [0, +\infty) \to [0,+\infty]$ such that 
$\psi(z,c), \ c > 0,$  has exactly two (counting multiplicity) negative zeros (say, $\lambda_1 \leq \lambda_0$)
in the half plane $\{\Re z < 0\}$ if and only if $c \in (0, C^*(\tau)]$.  Moreover, 
$C^*(\tau) = +\infty$ if and only if $\tau \leq 1/e$ while on the interval $(1/e,+\infty)$ function $C^*(\tau)$ is decreasing and $C^*(+\infty) =0$. Furthermore, $\Re \lambda_j < \lambda_1$ for every complex root of $\psi(\lambda_j,c)=0$.  Note also that $C^*(\tau)=c^*(\tau)$ for $\tau \in [0,\tau_1]$. 
\end{lem}
On the other hand, Lemma 17 (2-3), Remarks 19,20 in \cite{TT} and the proof of Lemma 10 in \cite{GT} imply 
\begin{lem} \label{PL} Let  $c > C^*(\tau)$.  Then every root $\lambda_j(c)$ of $\psi(z,c)=0$ is simple and  depends 
smoothly on $c$. Moreover, each vertical half-line $\Re z = a, \Im z \geq 0,$ contains at most one root $\lambda_j$ and all roots $\lambda_j, \Im \lambda_j \geq 0,$ can be ordered in such a way that  \quad 
$\dots < \Re \lambda_{j+1}(c) < \Re \lambda_{j}(c) < \dots < \Re \lambda_{0}(c) < \lambda_{-1}(c).$  Finally, if 
$\Re \lambda_j(c) \leq 0, \Im \lambda_j(c) \geq 0, $ then \  $c\tau\Im \lambda_j(c) \in (2j\pi,(2j+1)\pi)$. 
\end{lem}

\vspace{-2mm}

Here we would like to stress the following important fact: if $\lambda(c_0) =i w, w >0,$ is purely imaginary 
zero of $\psi(z,c_0)$ then $\exp(-ic\tau w) = -w^2-ic w$ and thus

\vspace{-5mm}

$$
\Re \lambda'(c_0) = \frac{2w^2(1+\tau w^2)}{c^2(1+\tau w^2)^2 + w^2(c^2\tau-2)^2} >0.
$$

\vspace{-3mm}

Therefore the point $\lambda(c) \in \C$ must  cross transversally the imaginary axis   {\it only from the left to the right}  and at the unique  moment  $c_0$.  In view of the above lemmas, this means that $\psi(z,c)$ can have more that one zero in 
$\Re z \geq 0$ if and only if at some {\it uniquely determined} moment $c=C^\star(\tau) > C^*(\tau)$ the point $c\tau\lambda_0(c)$  crosses the vertical segment $i[0,\pi]$.  Moreover, for each $c \leq C^\star(\tau)$ characteristic function $\psi(z,c)$ has 
only one zero (i.e. $\lambda_{-1}(c)$) in the open right half-plane, while for $c > C^\star(\tau)$ it has at least 
three zeros (i.e. $\lambda_{-1}(c), \ \lambda_{0}(c), \bar \lambda_{0}(c)$) in $\{\Re z > 0\}$. By the last lemma, 
the strip $\Pi:= (-\infty,0]\times [0,2\pi/(c\tau)] \subset \C$ contains at most one complex zero (i.e. $\lambda_0(c)$) for $c > C^*(\tau)$ and therefore $\Pi$ does not contain any zero of $\psi(z,c)$ for $c > C^\star(\tau)$. 

Hence, $C^\star$ can be determined as a unique positive real number such that equation $\psi(z,C^\star)=0$, or, equivalently, 
\begin{equation} \label{ct}
\cos (C^\star\tau w) = - w^2, \quad \sin (C^\star\tau w) = C^\star w. 
\end{equation}

\vspace{-4mm}

has a solution $z = iw$ with $C^\star\tau w \in [0,\pi]$.  
From the first equation of (\ref{ct}) we obtain that actually $C^\star\tau w \in (\pi/2,\pi]$. Therefore $1/\tau  = \sin(C^\star\tau w)/(C^\star\tau w) < 2/\pi$. This means that $C^\star(\tau)=+\infty$ for all $\tau \in [0,\pi/2]$.  On the other hand, if $\tau >\pi/2$, equation $1/\tau  = \sin(c\tau w)/(c\tau w)$ has a unique root $c\tau w$ on $(\pi/2,\pi]$ and therefore $w$ can be determined  uniquely as $\sqrt{-\cos(c\tau w)}$.  It is clear that also $w^4+c^2w^2=1$, from which $w^2(c) = 0.5(-c^2+ \sqrt{c^4+4})$.  This proves the representation (\ref{fot}).  Finally, it is easy to see that $\tau(c)$ strictly decreases on $(0, +\infty)$, with $\tau(+\infty)=\pi/2$. Therefore 
$$\tau > \tau(2) =\frac{\arccos(2-\sqrt{5})}{2\sqrt{\sqrt{5}-2}}=: \tau_2$$ implies that $\psi(z,c)$ with $c \geq 2$ has at least three zeros on $\Re z \geq 0$.

\vspace{-8mm}

\section*{Acknowledgments}

\vspace{-8mm}


This research was realized within the framework of the OPVK program, project CZ.1.07/2.300/20.0002. 
Sergei Trofimchuk was also partially supported  by FONDECYT (Chile), project 1071053, and by
CONICYT (Chile) through PBCT program ACT-56.

\end{document}